\documentclass[oneside]{amsart}
\usepackage{amsmath,ifthen,amscd,amssymb,
graphicx,subfigure}
\newcommand{\NSgeomfinite}{MR96c:20066}
\newcommand{\Epatterns}{MR1974064}
\newcommand{\Wise}{MR96m:20058}

\newcommand{\Cannon}{MR88a:20049}
\newcommand{\Thiel}{MR95e:20052}

\newcommand{\GS}{MR92c:20058}
\newcommand{\LS}{MR58:28182}
\newcommand{\Efftpft}{MR1950887}
\newcommand{\DShap}{MR92d:20055}
\newcommand{\ShapS}{MR96c:20068}
\graphicspath{{.}{./eps/}}

\newtheorem{thm}{Theorem}[section]

\newtheorem{propn}[thm]{Proposition}
\newtheorem{cor}[thm]{Corollary} 

\theoremstyle{definition}
\newtheorem{defn}[thm]{Definition}

\newtheorem{example}{Example}

\newcommand{\bi}{\begin{itemize}}
\newcommand{\ei}{\end{itemize}}
\newcommand{\be}{\begin{enumerate}}
\newcommand{\ee}{\end{enumerate}}
\newcommand{\bc}{\begin{center}}
\newcommand{\ec}{\end{center}}
\newcommand{\bt}{\begin{tabular}}
\newcommand{\et}{\end{tabular}}

\newcommand{\smallcaps}[1]{\textrm{\textsc{#1}}}

\newcommand{\cat}{\smallcaps{CAT}}

\newcommand{\ac}{almost convex}
\newcommand{\fftp}{falsification by fellow traveler property}
\newcommand{\cg}{Cayley graph}
\newcommand{\hnn}{HNN extension}

\newcommand{\gset}{generating set}

\begin{document}

\title[A non-Hopfian almost convex group]
  {A non-Hopfian almost convex group}

\author[M.~Elder]{Murray J. Elder}
\address{(As at Dec 2006) Dept of Mathematical Sciences\\
        Stevens Institute of Technology\\
        Hoboken NJ 07030 USA}
\email{melder@stevens.edu}

\date{Published in Journal of Algebra Vol 271, No. 1 (2004) 11--21}

\begin{abstract}
In this article we prove that an isometric multiple \hnn\ of a group
satisfying the \fftp\ is \ac. As a corollary, Wise's example of a
\cat(0) non-Hopfian group is \ac.
\end{abstract}

\maketitle

\section{Introduction}

In this article we give a proof that certain multiple \hnn s enjoying
certain geodesic conditions are \ac.  A non-Hopfian example of Wise
\cite{\Wise} satisfies these conditions, so we have a non-Hopfian
(hence non-residually finite) \ac\ group.

The article is organized as follows. In Section \ref{sec:def} we
define the properties of \ac ity and the \fftp, and state some facts
about them. We then define a multiple \hnn\, and the properties of
being strip equidistant and totally geodesic. A multiple \hnn\
satisfying both conditions and having finitely generated free
associated subgroups is said to be an {\em isometric multiple \hnn}.
In Section \ref{sec:egs} we give two simple examples of isometric
multiple \hnn s, including Wise's non-Hopfian example.  In Section
\ref{sec:thm} we state and prove the main theorem about \ac ity for
isometric multiple \hnn s.

The author wishes to thank Walter Neumann, Jon McCammond and Mike
Shapiro for their help, and the article's referee for her/his careful
reading and useful comments.

\section{Definitions}\label{sec:def}

Let $G$ be a group with finite \gset\ $X$, and let $\Gamma(G,X)$ be the corresponding \cg.
\begin{defn}[Almost convex]
$(G,X)$ is {\em \ac} if there is a constant $C>0$ such that for every pair of elements $g,g'$
 in the metric sphere of radius $N$ 
 at most distance 2 apart in $\Gamma(G,X)$ there is a path of length at most $C$ from $g$ to $g'$
which runs inside the metric ball of radius $N$.
\end{defn}
Cannon introduced the notion of an \ac\ group in \cite{\Cannon}, where he proved that
if a pair $(G,X)$ is \ac\ then there  is an efficient algorithm to construct any finite portion of its \cg.
Word hyperbolic groups, Coxeter groups and the fundamental groups of closed 3-manifolds with one of the eight geometries except Solvgeometry are \ac\ \cite{\DShap, \ShapS}.

Neumann and Shapiro gave a definition which they attributed to Cannon,
of a nice property that turns out to be relevant to \ac\ groups.
\begin{defn}[The \fftp]
  $(G,X)$ has the (asynchronous) {\em \fftp} if there is a constant
  $k>0$ such that every non-geodesic word is (asynchronously)
  $k$-fellow traveled by a shorter word in $\Gamma(G,X)$.
\end{defn}

\noindent
If a pair $(G,X)$ enjoys the asynchronous \fftp\ then it also enjoys
the synchronous \fftp\ \cite{\Efftpft}.  If $(G,X)$ has the \fftp\ 
then the full language of geodesics on $X$ is regular
\cite{\NSgeomfinite}.
 
\begin{propn}[Falsification by fellow traveling implies \ac]
\label{fftpimpliesac}
If $(G,X)$ has the \fftp\ then $(G,X)$ is \ac.
\end{propn}

\noindent
The proof of this can be found in \cite{\Efftpft}.  If $k$ is the
\fftp\ constant then the \ac\ constant $C$ is at most $3k$.  The
converse of Proposition \ref{fftpimpliesac} is false; in the present
article we prove that Wise's example is \ac, and in \cite{\Epatterns}
we prove that the full language of geodesics is not regular for the
same \gset, hence it doesn't enjoy the \fftp.

The enjoyment of either property is dependent on the choice of \gset\ 
\cite{\NSgeomfinite, \Thiel}.  The following result is proved in
\cite{\NSgeomfinite}, where the authors go on to prove that any
virtually abelian group has a \gset\ for which it enjoys the \fftp.
\begin{propn}[Abelian implies the \fftp]
\label{prop:abelian}
Any finite 
generating set for an abelian group has the \fftp.
\end{propn}

\begin{defn}[Multiple \hnn]
  Let $A$ be a group with finite \gset\ $X$ and relators $R$. Define
  an isomorphism $\phi_i:U_i \rightarrow V_i$ between pairs of
  isomorphic subgroups $U_i,V_i$, for $i\in [1,n]$. The {\em multiple
    \hnn} of $(A,X)$ with these isomorphisms is the group $G$ with
  presentation
$$
\langle X,s_1,\ldots , s_n|R, s_i^{-1}u_is_i=\phi_i(u_i); u_i\in
U_i, i\in[1,n] \rangle.$$
\end{defn}
The generators $s_i$ are called {\em stable letters}, and a subword of
the form $ s_i^{-1}u_is_i$ or $ s_iv_is_i^{-1}$ is called a {\em
  pinch}, where $u_i\in U_i, v_i\in V_i$.  A word that contains no
pinches is called {\em stable letter reduced}.  Britton's Lemma states
that if a freely reduced word containing stable letters is non-trivial
and represents the identity in the group then it must contain a pinch.
Two words are said to have {\em parallel stable letter structure} if
they have the exact same sequence of stable letters (when we ignore
the elements of the base group).  See \cite{\LS} for more details
about \hnn s.

When each $U_i$ is finitely generated, that is,
$U_i=\{u_{ij}:j=1,\ldots, m_i\}$, then define $
v_{ij}=\phi_i(u_{ij})$, and since $\phi_i$ is an isomorphism $V_i$ is
generated by $\{v_{ij}:j=1,\ldots, m_i\}$, and the multiple \hnn\ has
the presentation
$$\langle X,s_1, \ldots , s_n|R, s_i^{-1}u_{ij}s_i=v_{ij}; i\in[1,n],
j\in[1,m_i] \rangle,$$
which is finite when $R$ is finite.

If $X$ is an alphabet let $X^*$ denote the set of all words in the
letters of $X$ (including the empty word).
\begin{defn}[Geodesic, totally geodesic, strip equidistant]
  We say associated subgroups are {\em geodesic} if each freely
  reduced word in $\{u_{i1}^{\pm 1},\ldots, u_{i{m_i}}^{\pm 1}\}^*$
  and $\{v_{i1}^{\pm 1},\ldots, v_{i{m_i}}^{\pm 1}\}^*$ is geodesic,
  and {\em totally geodesic} if for each geodesic word $w\in U_i$
  [respectively $V_i$], $w \in \{u_{i1}^{\pm 1},\ldots, u_{im_i}^{\pm
    1}\}^*$ [respectively $w \in \{v_{i1}^{\pm 1},\ldots,
  v_{im_i}^{\pm 1}\}^*$].  Note that totally geodesic subgroups are
  geodesic.  We say the geodesic associated subgroups are {\em strip
    equidistant} if $|u_{ij}|=|v_{ij}|$ for each $i,j$.  Finally, we
  say a presentation for a multiple \hnn\ is (totally) geodesic
  [respectively strip equidistant] if all associated subgroups are.
\end{defn}

\begin{defn}[Isometric multiple \hnn]
  Let $A=\langle X|R\rangle$ be a finitely presented group with pairs
  of isomorphic finitely generated free subgroups $U_i,V_i$ for $i\in
  [1,n]$.  If the multiple \hnn\ of $(A,X)$ associating these
  subgroups is strip equidistant and totally geodesic, then we call it
  an {\em isometric multiple \hnn}.
\end{defn}

The \cg\ of such a presentation consists of copies of the \cg\ of the
base group $(A,X)$, glued together along the subspaces corresponding
to the free subgroups $U_i$ and $V_i$ by stable letter ``strips''. See
Figure \ref{fig:strip} in the next section for an illustration.

\section{Examples}\label{sec:egs}

In this section we give two examples of groups in the class of
isometric multiple \hnn s.
\begin{example}[Wise]\label{eg:wise}
  Let $A=\mathbb Z^2 \cong \langle a,b,c,d | c=ab,c=ba, d=c^2\rangle$
  and define the multiple \hnn\ $G_W$ by associating pairs of cyclic
  subgroups $\langle a\rangle, \langle d\rangle$ and $\langle
  b\rangle, \langle d\rangle$.  $G_W$ has the presentation
  $$\langle a,b,c,d,s,t|c=ab, c=ba,d=c^2, s^{-1}as=d,
  t^{-1}bt=d\rangle.$$
  Wise showed that this group is non-Hopfian and
  \cat(0) \raisebox{.8ex}{\hspace{.05mm}\normalfont \scriptsize 1}.
  \footnotetext[1]{Wise claimed in \cite{\Wise} that this group is
    automatic; the proof given was incorrect and its automaticity is
    as yet unresolved.}
  
  The \cg\ for $(G_W,\{a,b,c,d,s,t\})$ is the universal cover of the
  $1$-skeleton of the $2$-complex shown Figure \ref{fig:presW}. It is
  easy to check that this $2$-complex, metrized so that the edge
  labeled $c$ has half the length of the other edges, satisfies the
  link condition so is \cat(0).  This is the presentation that Wise
  used, minus the relation $d=c^2$, to prove \cat(0). Note that with
  this metric, the triangle corresponding to this relation is
  degenerate, so looks a little strange in the figure.
\begin{figure}
\bt{c}\includegraphics[scale=.3]{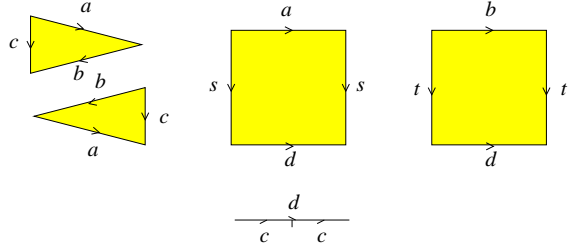}\et
\caption{The presentation $2$-complex for $(G_W,\{a,b,c,s,t\})$.\label{fig:presW}}
\end{figure}
The \cg\ can be viewed as being made up of ``planes'', corresponding
to copies of the \cg\ of the base group $(A,X)$, in this case
$(\mathbb Z^2,\{a,b,c,d\})$, glued together by ``strips'' made up of
copies of the squares (more generally metric rectangles) in Figure
\ref{fig:presW}, shown in Figure \ref{fig:strip}.
\begin{figure}[ht!]
  \begin{center}
      \includegraphics[width=9cm]{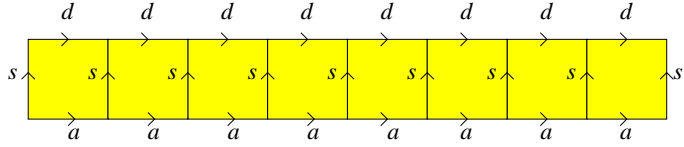}
  \end{center}
  \caption{A strip in $\Gamma (G_W,\{a,b,c,d,s,t\})$}
  \label{fig:strip}
\end{figure}

The presentation is by design strip equidistant, and totally geodesic,
since the words $a^n,b^n$ and $d^n$ are unique geodesic
representatives in $(\mathbb Z^2,\{a,b,c,d\})$ for elements of the
associated subgroups $\langle a\rangle, \langle b\rangle$ and $\langle
d\rangle$ respectively.

It follows that if a path in the \cg\ is not stable letter reduced,
then it can be shortened and therefore is not geodesic. Moreover two
geodesics from the identity to the same plane have parallel stable
letter structure, using Britton's Lemma.  In \cite{\Epatterns} we
consider in detail the geodesic structure of this presentation.
\end{example}

\begin{example}
Let $A=\mathbb F_2 \cong \langle a,b, | - \rangle$ and define the
(single) \hnn\ $G_2$ by defining an isomorphism between free subgroups
$ \langle a^2,b^3 | - \rangle$ and $ \langle b^2,aba | - \rangle$ by
$\phi(a^2)=b^2,\phi(b^3)=aba$.  $G_2$ has the presentation
$$\langle a,b,s| s^{-1}a^2s=b^2,
s^{-1}b^3s=aba\rangle.$$

The presentation $2$-complex shown in Figure \ref{fig:brick} is easily
seen to satisfy the link condition, so is \cat(0). Moreover it can be
viewed as a \cat(0) squared complex, so by Gersten-Short \cite{\GS} is
biautomatic.
\begin{figure}
\bt{c}\includegraphics[scale=1]{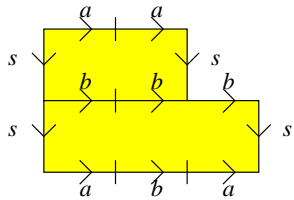}\et
\caption{The presentation $2$-complex for $G_2$.\label{fig:brick}}
\end{figure}

It is easily checked that the associated subgroups are totally
geodesic, and that the presentation is strip equidistant, so
$(G_2,\{a,b,s\})$ is an isometric (multiple) \hnn.
\end{example}

\section{The Main theorem}\label{sec:thm}

We now prove that an  isometric multiple \hnn\ is \ac.
\begin{thm}
Let $(G, X\cup \{s_i\}_{i=1} ^n)$ be an isometric multiple \hnn\ with
base group $(A,X)$.  If $(A,X)$ has the \fftp\ then $G$ is \ac\ with
respect to the \gset\ $X\cup \{s_i\}_{i=1} ^n$.
\end{thm}

\begin{proof}
Let $S(N)$ denote the metric sphere of radius $N$ and $B(N)$ the
metric ball of radius $N$ in $\Gamma(G, X\cup \{s_i\}_{i=1} ^n)$.  Let
$g,g'\in S(N)$ with $d(g,g')\leq 2$ realized by a path $\gamma$.  Let
$w,w'$ be geodesic words for $g,g'$ respectively.  Since the
presentation is strip equidistant, $w$ and $w'$ are stable letter reduced.

Note that the metric sphere [ball] of radius $N$ in the base group
$(A,X)$ is a subset of the sphere [ball] of radius $N$ in $\Gamma(G,
X\cup \{s_i\}_{i=1} ^n)$.  Without loss of generality we may assume
the \fftp\ constant $k$ for $(A,X)$ is an even integer. We may also
assume that $X$ is inverse closed. Denote the endpoint of a path $w$
from the identity by $\overline w$.

The argument is divided into 3 cases.

\subsection*{Case 1} 
$w\gamma (w')^{-1}$ has no stable letters. Then since $(A, X)$ has the
\fftp\, we are done. (Recall that $C=3k$.)

\subsection*{Case 2}
$\gamma$ involves a stable letter. That is, $\gamma =s,sx$ or $xs$ for
$s\in \{s_i^{\pm 1}\}_{i=1}^n$ and $x$ any generator (or inverse of a
generator) except $s^{-1}$.  By Britton's Lemma $w\gamma (w')^{-1}$
contains a pinch so either $w$ or $w'$ has an $s^{-1}$.  Without loss
of generality assume $w=w_1s^{-1}w_2$ where $s^{-1}w_2s$ is a pinch.

If $\gamma =s$ or $\gamma= sx$ then the path $w_1s^{-1}w_2$ is shown
in Figure \ref{fig:sws}(a).
\begin{figure}[ht!]
  \begin{center}
         \subfigure[$\gamma =s, sx$]{
         \includegraphics[width=4cm]{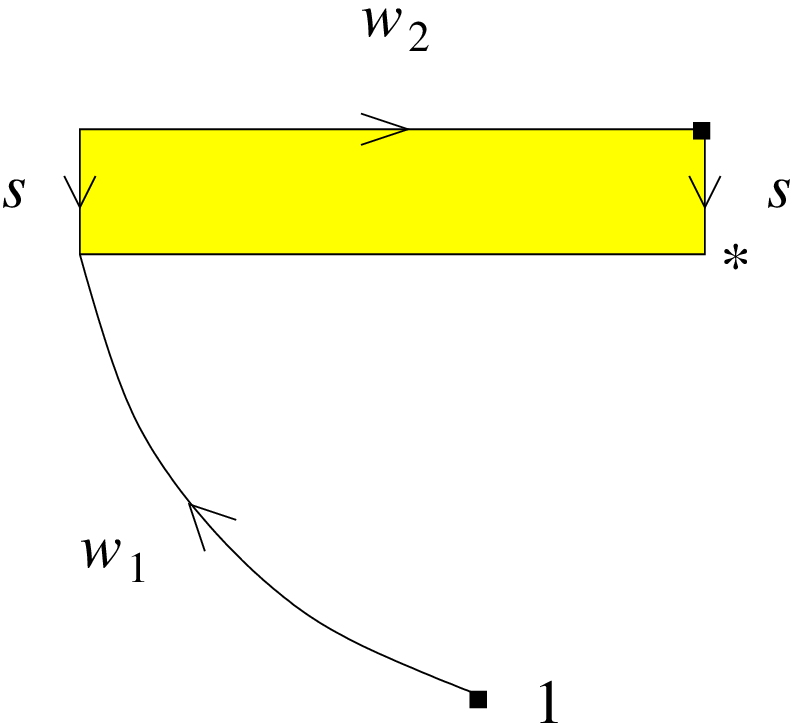}}
        \qquad \qquad
       \subfigure[$\gamma =xs$]{
 \includegraphics[width=4cm]{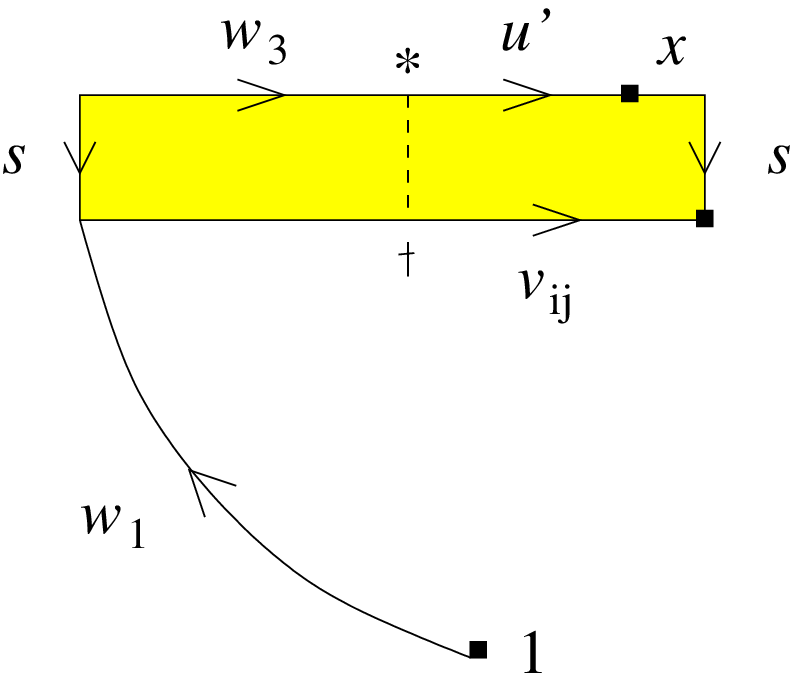}}
        \qquad \qquad
       \subfigure[$\gamma =xs$]{
      \includegraphics[width=4cm]{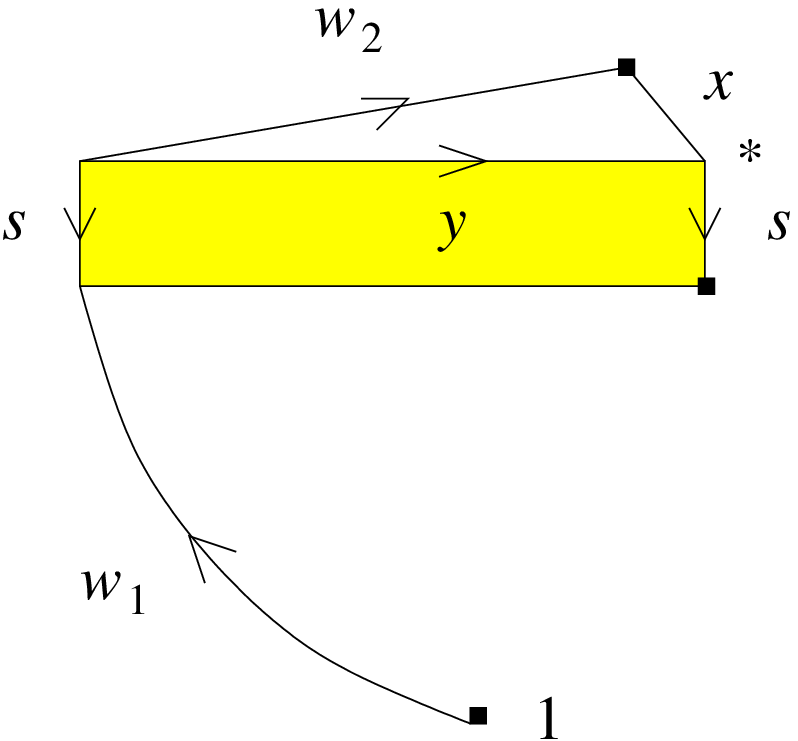}}
        \qquad \qquad
    \subfigure[$\gamma =xs : |w_2|=|y|$]{
      \includegraphics[width=4cm]{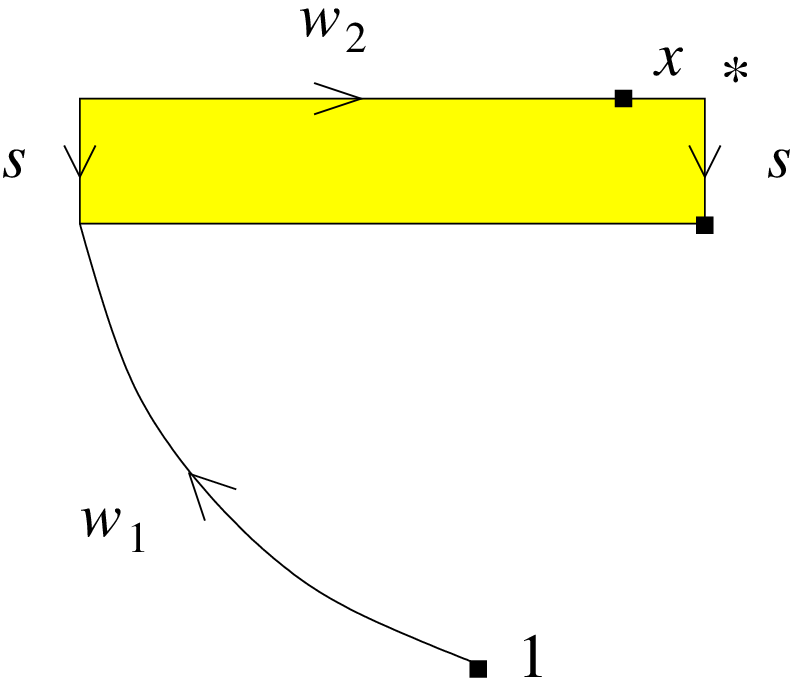}}
   \end{center}
  \caption{Case 2}
 \label{fig:sws} 
\end{figure}
The subword $w_2$ is geodesic, and since $s^{-1}w_2s$ is a pinch, it
is an element of an associated subgroup $U_i$ or $V_i$.  Since the
presentation is totally geodesic, this path must run along the top of
a strip, and the bottom of the strip is a word of the same length.
Thus the point labeled $\ast$ (which is $\overline {ws}$) lies in
$B(N-1)$, so if $\gamma=sx$ then $\gamma$ lies in $B(N)$ and if
$\gamma=s$ we have a contradiction.
 
 If $\gamma=xs$ then $w_2x$ evaluates to an element of an associated
 subgroup.  If $x$ is a stable letter
 then we are in the previous case.  So we may assume $w_2x\in X^*$.
 If $w_2x$ is geodesic then it must run along the strip, so without
 loss of generality there is some $u_{ij}$ ending in $x$, with
 $s^{-1}u_{ij}s=v_{ij}$. Let $u_{ij}=u'x$ and $w_2=w_3u'$. We show the
 path $w_1s^{-1}w_3u'$ in Figure \ref{fig:sws}(b), and draw the
 subword $w_2x=w_3u_{ij}=w_3u'x$ on top of a strip.  Consider the path
 $(u')^{-1}sv_{ij}$ between $g$ and $g'$. The point labeled $\ast$
 lies on the path $w$ at distance $N-|u'|$ from the identity.  The
 point labeled $\dag$ is distance $N-|w_2|-1+|w_3|=N-|u'|-1$ from the
 identity, since it can be reached by traveling along $w_1$ then along
 the bottom of the strip shown in Figure \ref{fig:sws}(b). Thus the
 path $(u')^{-1}sv_{ij}$ between $g$ and $g'$ stays inside $B(N)$. It
 has length at most $2\max\{|u_{ij}|:u_{ij}$ is a generator of $U_i,
 i\in [1,n]\}$.
 
 If $w_2x$ is not geodesic then let $y$ be the geodesic for $w_2x$ as
 in Figure \ref{fig:sws}(c).  Now $|y|<|w_2x|=|w_2|+1$ so $|y|\leq
 |w_2|$.  There is a path $w_1sy$ to the point labeled $\ast$ in
 Figure \ref{fig:sws}(c) of length $|w_1|+1+|y|\leq |w_1|+1+|w_2|=N$
 so if the point $\ast \in S(N)$ we have Case 2 (with $\gamma = s$, if
 $\ast \in S(N-1)$ then $\gamma$ lies in $ B(N)$ and if $\ast \in
 B(N-2)$ then we have a contradiction.

\subsection*{Case 3}
$\gamma$ has no stable letters and $w$ has a stable letter.  Then
$w,w'$ have parallel stable letter structure.  Suppose $s\in
\{s_i^{\pm 1}\}^n_{i=1}$ is the last stable letter of $w$, so
$w=w_1sw_2, w'=w'_1sw'_2$, where $w_2$ and $w_2'$ have no stable
letters.

We will draw the portion of the \cg\ which contains the paths
$w_2,\gamma$ and $w_2'$ in Figures \ref{fig:uy} and \ref{fig:w2gamu2}.
The path $w_2\gamma(w'_2)^{-1}$ starting at the point labeled $\ast$
in both figures is a word in $X^*$, and moreover evaluates to an
element of an associated subgroup. Therefore it evaluates to a word
that runs along the top of a strip.  This strip is shown in each
figure, and below it there are two paths (not shown) $w_1$ and $w_1'$
back to the identity.

If $w_2\gamma(w'_2)^{-1}$ is geodesic then it runs along the top of
the strip. The path $\gamma$ could lie on top of at most two
rectangles in the strip, so by tracing around these rectangles we stay
within $B(N)$. Therefore we can find a path inside $B(N)$ from $g$ to
$g'$ of length at most $4\max\{|u_{ij}|:u_{ij}$ is a generator of
$U_i, i\in [1,n]\}$.

If $w_2\gamma(w'_2)^{-1}$ is not geodesic then by the \fftp\ in
$(A,X)$ there is a shorter word $u$ which synchronously $k$-fellow
travels it. If $u$ is geodesic, it runs along the strip, and we show
these paths in Figure \ref{fig:uy}.
\begin{figure}[ht!]
  \begin{center}
         \includegraphics[width=8cm]{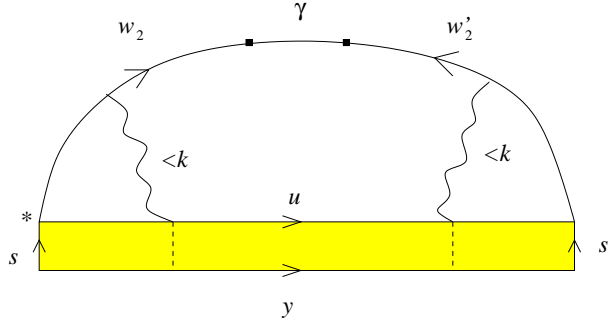}
  \end{center}
  \caption{Case 3: $u$ is geodesic.}
  \label{fig:uy}
\end{figure}
Let $y=sus^{-1}$ be the geodesic on the other side of the strip from
$u$.  Then $y$ has length at most $|w_2\gamma(w'_2)^{-1}|-1$, starts
at $\overline {w_1} \in S(N-|w_2|-1)$ and ends at $\overline {w_1'}
\in S(N-|w_2'|-1)$, so lies in $B(N)$.

Let $w(t)$ denote the point at distance $t\in \mathbb R_{\geq 0}$
along the path $w$ from its start point, where $w(t)=\overline w$ for
$t\geq |w|$.

If $|w_2|,|w_2'| \geq \frac{k}{2}$ then consider the path starting at
$\overline w=w(N)$ retracing along $w$ to $w(N-\frac{k}{2})$.  The
paths $w_2\gamma (w_2')^{-1}$ and $u$ based at $\ast$ in Figure
\ref{fig:uy} $k$-fellow travel in the base group $(A,X)$, so the paths
$w_1sw_2\gamma (w_2')^{-1}$ and $w_1su$ (based at the identity)
$k$-fellow travel is the larger group.  Thus there is a path of length
at most $k$ from $w(N-\frac{k}{2})$ to $w_1su(N-\frac{k}{2})$ which
lies in $B(N)$. This path is drawn as a jagged line in the figure.
From here cross to $w_1y(N-\frac{k}{2})$ by an edge $s^{-1}$, and
travel along $y$ (which is inside $B(N)$) to $w_1y(N+|\gamma|+
\frac{k}{2})$, then cross the strip by an edge $s$ to
$w_1su(N+|\gamma|+ \frac{k}{2})$.  Again since $w_1sw_2\gamma
(w_2')^{-1}$ and $w_1su$ $k$-fellow travel there is a path of length
at most $k$ to $w_1sw_2\gamma (w_2')^{-1}(N+|\gamma|+ \frac{k}{2})$.
Now $w_1sw_2\gamma (w_2')^{-1}(N+|\gamma|+
\frac{k}{2})=w'(N-\frac{k}{2}) \in S(N-\frac{k}{2})$ and
$w_1su(N+|\gamma|+ \frac{k}{2})=w_1'su^{-1}(N- \frac{k}{2}-c)$ for
some positive $c$ since $u$ is shorter than $w_2\gamma (w_2')^{-1}$,
so the path between these points lies in $B(N)$.  Finally, travel
along $w'$ to $\overline{w'}$.  The concatenation of these paths stays
within $B(N)$, and has length at most
$\frac{k}{2}+k+(2\frac{k}{2}+2)+k+\frac{k}{2}=4k+4$.

If $|w_2|,|w_2'| < \frac{k}{2}$ then we have a path from $\overline w$
to $\overline {w_1}$, along $y$ to $\overline {w_1'}$, then along $w'$
to $\overline {w'}$, which runs inside $B(N)$ and has length less than
$\frac{k}{2}+1+(2\frac{k}{2}+2)+1+\frac{k}{2}=2k+4.$

If $|w_2| < \frac{k}{2}, |w_2'| \geq \frac{k}{2}$ (or vice versa) then
a combination of the above arguments gives a path inside $B(N)$ of
length at most $\frac{k}{2}+1+(2\frac{k}{2}+2)+1+k+\frac{k}{2}=3k+4$.
Namely, travel back along $w$ to $\overline{w_1}$, then along $y$ to
$w_1y(N+|\gamma|+\frac{k}{2})$, across the strip to
$w_1u(N+|\gamma|+\frac{k}{2})$ and by a path of length at most $k$ to
$w_1sw_2\gamma (w_2')^{-1}(N+|\gamma|+
\frac{k}{2})=w'(N-\frac{k}{2})$, then along $w'$ to $\overline{w'}$.
The other case is similar.

 If $u$ is not geodesic then by the \fftp\ there is a shorter word $v$
 in the base group $(A,X)$ which synchronously $k$-fellow travels $u$
 as shown in Figure \ref{fig:w2gamu2}.
\begin{figure}[ht!]
  \begin{center}
         \includegraphics[width=8cm]{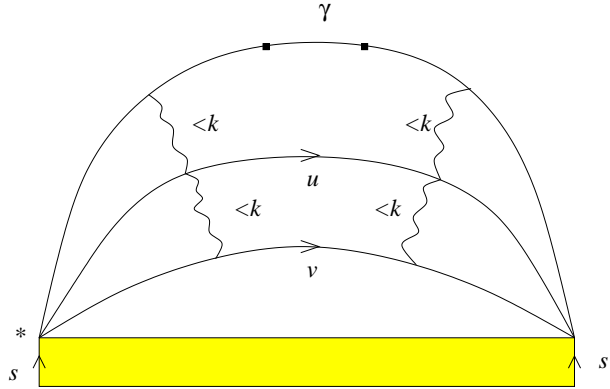}
  \end{center}
  \caption{Case 3: $u$ is not geodesic.}
  \label{fig:w2gamu2}
\end{figure}
Then $|v|\leq|w_2\gamma (w_2')^{-1}|-2\leq |w_2|+|w_2'|$ so $v$ lies
in $B(N)$.

If $|w_2|,|w_2'| \geq \frac{k}{2}$ then there is a path starting at
$\overline{w}=w(N)$, back along $w$ to $w(N-\frac{k}{2})$. From here
there is a path of length at most $k$ to $w_1su(N-\frac{k}{2})$ then
from here a path of length at most $k$ to $w_1sv(N-\frac{k}{2})$.
These paths lie in $B(N)$. From here travel along $v$ (which lies in
$B(N)$) to $w_1sv(N+|\gamma|+\frac{k}{2})$. From here there is a path
of length at most $k$ to $w_1su(N+|\gamma|+\frac{k}{2})$ and another
to $w_1sw_2\gamma (w_2')^{-1}(N+|\gamma|+\frac{k}{2})$.  The point
$w_1sv(N+|\gamma|+\frac{k}{2})=w_1'sv^{-1}(N-\frac{k}{2}-c)$ for some
positive $c$ and
$w_1su(N+|\gamma|+\frac{k}{2})=w_1'su^{-1}(N-\frac{k}{2}-d)$ for some
positive $d$, so these paths lie in $B(N)$.  Finally travel from
$w_1sw_2\gamma (w_2')^{-1}(N+|\gamma|+\frac{k}{2})=w'(N-\frac{k}{2})$
back to $\overline{w'}$.  The total length of the entire path is at
most $\frac{k}{2}+k+k+(2\frac{k}{2}+2) +k+k +\frac{k}{2}=6k+2$.

If $|w_2|,|w_2'| < \frac{k}{2}$ then $|v|\leq |w_2|+|w_2'|<k$ so there
is a path from $\overline w$ to $\overline {w_1}$, then along $v$ to
$\overline {w_1's}$ to $\overline {w's}$, which runs inside $B(N)$ and
has length less than $2k$.

If $|w_2| < \frac{k}{2}, |w_2'| \geq \frac{k}{2}$ (or vice versa) then
a combination of the previous two arguments gives a path from $
\overline w$ to $\overline{w_1s}$, then along $v$ to
$w_1v(N+|\gamma|+\frac{k}{2})$. From here there is a path to
$w_1u(N+|\gamma|+\frac{k}{2})$ of length at most $k$ which lies in
$B(N)$, then a path of length at most $k$ to $w_1sw_2\gamma
(w_2')^{-1}(N+|\gamma|+ \frac{k}{2})=w'(N-\frac{k}{2})$, then along
$w'$ to $\overline{w'}$, and its total length is at most
$\frac{k}{2}+(2\frac{k}{2}+2)+k+k+\frac{k}{2}=4k+2.$ The other case is
similar.

This completes all possible cases.  The \ac ity constant $C$ is at
most $\max\{6k+2, |u_{ij}|:u_{ij}$ is a generator of $U_i, i\in
[1,n]\}$ where $k$ is the \fftp\ constant for $(A,X)$.
\end{proof}

\begin{cor}
Wise's example  is \ac\ and non-Hopfian.
\end{cor}
\begin{proof}
  We showed in Example \ref{eg:wise} that Wise's example is an
  isometric multiple \hnn, and the base group is abelian so by
  Proposition \ref{prop:abelian} satisfies the \fftp. Wise proved that
  the example is non-Hopfian in \cite{\Wise}.
\end{proof}

\bibliography{refs}
\bibliographystyle{plain}

\end{document}